\documentclass{amsart}

\usepackage{fancyhdr}
\usepackage{lastpage}
\usepackage{stmaryrd,yhmath}

\newcommand{\bdis}{\begin{displaymath}}
\newcommand{\edis}{\end{displaymath}}
\newcommand{\be}{\begin{equation}}
\newcommand{\ee}{\end{equation}}
\newcommand{\mbb}{\mathbb}
\newcommand{\mcal}{\mathcal}

\newcommand{\vp}{\varphi}

\newcommand{\vth}{\vartheta}

\newcommand{\zf}{\zeta\left(\frac{1}{2}+it\right)}

\theoremstyle{definition}

\theoremstyle{remark}
\newtheorem{remark}[]{Remark}

\newtheorem*{mydef1}{{\bf Theorem}}

\newtheorem*{mydef4}{{\bf Corollary}}

\numberwithin{equation}{section}

\begin{document}

\title{Jacob's ladders and invariant set of constraints for the reversely iterated integrals (energies) in the theory
of the Riemann zeta-function}

\author{Jan Moser}

\address{Department of Mathematical Analysis and Numerical Mathematics, Comenius University, Mlynska Dolina M105, 842 48 Bratislava, SLOVAKIA}

\email{jan.mozer@fmph.uniba.sk}


\begin{abstract}
In this paper we obtain an extension of the set of non-local equalities by adding to it new set of local equalities.
Namely, we obtain an invariant set of equalities on the set of reversely iterated integrals (energies). In other
words, we obtain a new continuum set of constraints on behaviour of the function
$\zf,\ t\to\infty$.
\end{abstract}
\maketitle

\section{Introduction}

\subsection{}

Let us remind we have proved the following theorem (see \cite{4}): for the class $G(S)$ of functions
\bdis
g=g(u_1,\dots,u_n),\ (u_1,\dots,u_n)\in S\subset{R}^n
\edis
such that
\be \label{1.1}
g\geq 0,\ g=o\left(\frac{T}{\ln T}\right),\ T\to\infty
\ee
we have the following formula
\be \label{1.2}
\forall\ g\in G(S):\
\int_{\overset{k}{T}}^{\overset{k}{\wideparen{T+g}}}\prod_{r=0}^{k-1}\tilde{Z}^2[\vp_1^r(t)]{\rm d}t=g,\
k=1,\dots,k_0
\ee
for every fixed $k_0\in\mbb{N}$ and for every sufficiently big $T>0$.

\begin{remark}
Let us notice that the level of precision of the formula (\ref{1.2}) is characterized sufficiently by the following
example
\bdis
\int_{\overset{257}{T}}^{\overset{257}{\wideparen{T+10^{-60}}}}\prod_{r=0}^{256}
\left|\zeta\left(\frac 12+i\vp_1^r(t)\right)\right|^2{\rm d}t=
10^{-60}\left\{ 1+\mcal{O}\left(\frac{\ln\ln T}{\ln T}\right)\right\}\ln^{257}T.
\edis
This formula is not accessible neither by the Riemann-Siegel formula
\bdis
Z(t)=2\sum_{n\leq\bar{t}}\frac{1}{\sqrt{n}}\cos\{\vth(t)-t\ln n\}+\mcal{O}(t^{-1/4}),\
\bar{t}=\sqrt{\frac{t}{2\pi}},
\edis
nor by the current methods in the theory of the Riemann zeta-function.
\end{remark}

The formula (\ref{1.2}) appears within the following context.

\subsection{}

The complicated signal
\bdis
\prod_{r=0}^{k-1}\tilde{Z}^2[\vp_1^r(t)];\ \tilde{Z}^2[\vp_1^0(t)]=\tilde{Z}^2(t)
\edis
is generated by the primary signal
\be \label{1.3}
\begin{split}
 & Z(t)=e^{i\vth(t)}\zf , \\
 & \vth(t)=-\frac t2\ln\pi+\mbox{Im}\ln\Gamma\left(\frac 14+i\frac t2\right)
\end{split}
\ee
which itself is generated by the Riemann zeta-function on the critical line. Namely, in connection with
(\ref{1.3}), we have introduced (see \cite{2}, (9.1), (9.2)) the formula
\be \label{1.4}
\tilde{Z}^2(t)=\frac{{\rm d}\vp_1(t)}{{\rm d}t},
\ee
where
\be \label{1.5}
\begin{split}
 & \tilde{Z}^2(t)=\frac{Z^2(t)}{2\Phi'_\vp[\vp(t)]}=\frac{\left|\zf\right|^2}{\omega(t)}, \\
 & \omega(t)=\left\{ 1+\mcal{O}\left(\frac{\ln\ln t}{\ln t}\right)\right\}\ln t.
\end{split}
\ee
The function
\bdis
\vp_1(t)
\edis
which is called the Jacob's ladder (see our paper \cite{1}) according to the Jacob's dream in Chumash,
Bereishis, 28:12, has the following properties:
\begin{itemize}
 \item[(a)]
 \bdis
 \vp_1(t)=\frac 12\vp(t),
 \edis
 \item[(b)] the function $\vp(t)$ is solution of the non-linear integral equation (see \cite{1}, \cite{2})
 \bdis
 \int_0^{\mu[x(T)]}Z^2(t)e^{-\frac{2}{x(T)}t}{\rm d}t=\int_0^T Z^2(t){\rm d}t,
 \edis
 where each admissible function $\mu(y)$ generates the solution
 \bdis
 y=\vp_\mu(T)=\vp(T);\ \mu(y)\geq 7y\ln y.
 \edis
\end{itemize}

\begin{remark}
The main goal of introducing Jacob's ladders is described in \cite{1}, where we have proved by making use of these
Jacob's ladders, that the Hardy-Littlewood integral (1918)
\bdis
\int_0^T\left|\zf\right|^2{\rm d}t
\edis
has -- in addition to the Hardy-Littlewood expression (and also other similar to this one) possessing an
unbounded error at $T\to\infty$ -- the following set of almost exact expressions
\bdis
\begin{split}
& \int_0^T\left|\zf\right|^2{\rm d}t=\vp_1(T)\ln\vp_1(T)+(c-\ln 2\pi)\vp_1(T)+ \\
& + c_0+
\mcal{O}\left(\frac{\ln T}{T}\right),\ T\to\infty,
\end{split}
\edis
where $c$ is the Euler's constant, and $c_0$ is the constant from the Titchmarsh-Kober-Atkinson formula
(see \cite{5}, p. 141).
\end{remark}

\begin{remark}
 The Jacob's ladder $\vp_1(T)$ can be interpreted by our formula (see \cite{1}, (6.2))
 \be \label{1.6}
 T-\vp_1(T)\sim (1-c)\pi(T);\ \pi(T)\sim\frac{T}{\ln T},\ T\to\infty,
 \ee
 where $\pi(T)$ is the prime-counting function, as an asymptotically complementary function to the function
 \bdis
 (1-c)\pi(T)
 \edis
 in the following sense
 \be \label{1.7}
 \vp_1(T)+(1-c)\pi(T)\sim T,\ T\to\infty.
 \ee
\end{remark}

\begin{remark}
Let us notice explicitly that the main reason for the substitution
\bdis
\frac{T}{\ln T}\longrightarrow \pi(T),\ T\to\infty
\edis
by the Hadamard-de la Vall\' e Poussin formula
\bdis
\pi(T)\sim \frac{T}{\ln T},\ T\to\infty
\edis
lies in availability of the asymptotic law (\ref{1.7}) of complementarity of the functions
\bdis
\vp_1(T),\ (1-c)\pi(T),\ T\to\infty,
\edis
(see also \cite{1}, (6.1) -- (6.6)).
\end{remark}

\section{On structure of non-local equalities}

\subsection{}

Now we turn back to the formula (\ref{1.2}). It is clear that the following set of equalities
\be \label{2.1}
\begin{split}
 & \int_{\overset{1}{T}}^{\overset{1}{\wideparen{T+g}}}\tilde{Z}^2(t){\rm d}t=
 \int_{\overset{2}{T}}^{\overset{2}{\wideparen{T+g}}}\prod_{r=0}^1\tilde{Z}^2[\vp_1^r(t)]{\rm d}t= \dots =\\
 & = \int_{\overset{k}{T}}^{\overset{k}{\wideparen{T+g}}}\prod_{r=0}^{k-1}\tilde{Z}^2[\vp_1^r(t)]{\rm d}t,\
 k=1,2,\dots,k_0
\end{split}
\ee
is contained in (\ref{1.2}) for every sufficiently big $T$. The product
\be \label{2.2}
\prod_{r=0}^{k-1}\tilde{Z}^2[\vp_1^r(t)]
\ee
in (\ref{2.1}) has the following properties. First of all we have that (see \cite{3}, (2.10))
\be \label{2.3}
\begin{split}
 & \vp_1^0(t)=t\in [\overset{k}{T},\overset{k}{\wideparen{T+g}}] \ \Rightarrow \
 \vp_1^r(t)\in [\overset{k-r}{T},\overset{k-r}{\wideparen{T+g}}] , \\
 & r=1,\dots,k,\ k\leq k_0,
\end{split}
\ee
i.e.
\be \label{2.4}
\begin{split}
 & \vp_1^1(t)\in [\overset{k-1}{T},\overset{k-1}{\wideparen{T+g}}] , \\
 & \vp_1^2(t)\in [\overset{k-2}{T},\overset{k-2}{\wideparen{T+g}}] , \\
 & \vdots \\
 & \vp_1^{k-1}(t)\in [\overset{1}{T},\overset{1}{\wideparen{T+g}}], \\
 & \vp_1^k(t)\in [\overset{0}{T},\overset{0}{\wideparen{T+g}}]=[T,T+g].
\end{split}
\ee
Next, we have the following properties for the segments in (\ref{2.3}), (see \cite{3}, (2.5) -- (2.7))
\bdis
g=o\left(\frac{T}{\ln T}\right) \ \Rightarrow
\edis
\be \label{2.5}
|[\overset{k}{T},\overset{k}{\wideparen{T+g}}]|=\overset{k}{\wideparen{T+g}}-\overset{k}{T}=o\left(\frac{T}{\ln T}\right) ,
\ee
\be \label{2.6}
|[\overset{k-1}{T+g},\overset{k}{\wideparen{T}}]|=\sim (1-c)\pi(T);\ \pi(T)\sim\frac{T}{\ln T},
\ee
\be \label{2.7}
[T,T+g]\prec [\overset{1}{T},\overset{1}{\wideparen{T+g}}]\prec \dots \prec
[\overset{k}{T},\overset{k}{\wideparen{T+g}}],
\ee
where, of course, (see (\ref{2.6}))
\be \label{2.8}
\rho\{[\overset{r-1}{T},\overset{r-1}{\wideparen{T+g}}],[\overset{r}{T},\overset{r}{\wideparen{T+g}}]\}\sim
(1-c)\pi(T)
\ee
and $\rho$ is the distance of corresponding segments.

\begin{remark}
Asymptotic behaviour of the following disconnected set
\be \label{2.9}
\Delta(T,k,g)=\bigcup_{r=0}^k [\overset{r}{T},\overset{r}{\wideparen{T+g}}]
\ee
(see (\ref{2.5})--(\ref{2.7}), comp. \cite{3}, (2.9)) is as follows: if $T\to\infty$ then the components of the set
(\ref{2.9}) recede unboundedly each from other and all together are receding to infinity. Hence, if $T\to\infty$ then the
set (\ref{2.9}) looks like one dimensional Friedmann-Hubble expanding universe.
\end{remark}

\begin{remark}
Consequently, we notice the following:
\begin{itemize}
 \item[(a)] the equalities in (\ref{2.1}) are non-local ones (see (\ref{2.7}), (\ref{2.8})) -- this is the property
 of external non-localization,
 \item[(b)] simultaneously, the integrals in (\ref{1.2}) with $k=2,\dots,k_0$ are non-local (see (\ref{1.2}),
 (\ref{2.2})--(\ref{2.4})) -- this is the property of internal non-localization.
\end{itemize}
\end{remark}

In this paper we obtain an extension of the set of non-local equalities (\ref{2.1}) by adding to it a new set of
local equalities -- i.e. equalities without the property (a) -- for the reversely iterated integrals (energies). In other
words, we obtain an invariant set (at $T\to\infty$) of equalities-constraints on the set of reversely
iterated integrals (energies).

\section{Theorem on a set of local equalities}

\subsection{}

Let us remind (see \cite{3}, (5.2), (5.7)) that
\be \label{3.1}
T<\overset{1}{T}<\overset{2}{T}<\dots<\overset{k-1}{T},\ k=2,\dots,k_0,
\ee
where
\be \label{3.2}
\vp_1(\overset{k}{T})=\overset{k-1}{T}.
\ee
Next, putting
\bdis
T\longrightarrow \overset{1}{T}
\edis
in (\ref{3.1}) we obtain that
\be \label{3.3}
\overset{1}{T}<\overset{2}{T}<\dots<\overset{k}{T} .
\ee
Now, let us consider the segment (see (\ref{1.2}))
\bdis
[\overset{1}{T},\overset{1}{\wideparen{T+g}}]
\edis
and assign to it (for every fixed $g\in (0,\infty)$) the following sequence of segments
\be \label{3.4}
[\overset{2}{T},\overset{1}{\wideparen{\overset{1}{T}+g}}],\ [\overset{3}{T},\overset{1}{\wideparen{\overset{2}{T}+g}}],\dots ,
[\overset{k}{T},\overset{1}{\wideparen{\overset{k-1}{T}+g}}],
\ee
and, consequently, to each of these segments we assign the integral (energy)
\bdis
[\overset{p}{T},\overset{1}{\wideparen{\overset{p-1}{T}+g}}] \longrightarrow \int_{\overset{p}{T}}^{\overset{1}{\wideparen{\overset{p-1}{T}+g}}}
\tilde{Z}^2(t){\rm d}t,\ p=2,\dots,k.
\edis
Since
\bdis
\int_{\overset{1}{T}}^{\overset{1}{\wideparen{T+g}}}\tilde{Z}^2(t){\rm d}t=g
\edis
for every sufficiently big $T$, then (see (\ref{3.2}), (\ref{3.3}))
\be \label{3.5}
\begin{split}
& \int_{\overset{1}{T}}^{\overset{1}{\wideparen{T+g}}}\tilde{Z}^2(t){\rm d}t=
\int_{\overset{1}{\wideparen{\overset{p-1}{T}}}}^{\overset{1}{\wideparen{\overset{p-1}{T}+g}}}\tilde{Z}^2(t){\rm d}t= \\
& =\int_{\overset{p}{T}}^{\overset{1}{\wideparen{\overset{p-1}{T}+g}}}\tilde{Z}^2(t){\rm d}t=g,\ p=2,\dots,k.
\end{split}
\ee
Hence, we have the following formula (see (\ref{1.2}), (\ref{3.5}))
\be \label{3.6}
\int_{\overset{p}{T}}^{\overset{1}{\wideparen{\overset{p-1}{T}+g}}}\tilde{Z}^2(t){\rm d}t=
\int_{\overset{p}{T}}^{\overset{p}{\wideparen{T+g}}}\prod_{r=0}^{p-1}\tilde{Z}^2[\vp_1^r(t)]{\rm d}t=g,\  p=2,\dots,k.
\ee

\subsection{}

By the similar way we obtain
\be \label{3.7}
\begin{split}
& [\overset{p}{T},\overset{2}{\wideparen{\overset{p-2}{T}+g}}] \longrightarrow
\int_{\overset{p}{T}}^{\overset{2}{\wideparen{\overset{p-2}{T}+g}}}\prod_{r=0}^1=\int_{\overset{p}{T}}^{\overset{p}{\wideparen{T+g}}}\prod_{r=0}^{p-1}=g, \ p=3,\dots,k, \\
& [\overset{p}{T},\overset{3}{\wideparen{\overset{p-3}{T}+g}}] \longrightarrow \int_{\overset{p}{T}}^{\overset{3}{\wideparen{\overset{p-3}{T}+g}}}\prod_{r=0}^2=
\int_{\overset{p}{T}}^{\overset{p}{\wideparen{T+g}}}\prod_{r=0}^{p-1}=g,\ p=4,\dots,k, \\
& \vdots \\
& [\overset{p}{T},\overset{k-2}{\wideparen{\overset{p-k+2}{T}+g}}] \longrightarrow \int_{\overset{p}{T}}^{\overset{k-2}{\wideparen{\overset{p-k+2}{T}+g}}}\prod_{r=0}^{k-3}=
\int_{\overset{p}{T}}^{\overset{p}{\wideparen{T+g}}}\prod_{r=0}^{p-1}=g,\ p=k-1,k, \\
& [\overset{p}{T},\overset{k-1}{\wideparen{\overset{p-k+1}{T}+g}}] \longrightarrow
\int_{\overset{p}{T}}^{\overset{k-1}{\wideparen{\overset{p-k+1}{T}+g}}}\prod_{r=0}^{k-2}=
\int_{\overset{p}{T}}^{\overset{p}{\wideparen{T+g}}}\prod_{r=0}^{p-1}=g,\ p=k.
\end{split}
\ee

\begin{mydef1}
For every fixed (comp. (\ref{1.1}))
\bdis
g\in (0,+\infty),\ k_0\in\mbb{N},\ k_0\geq 2
\edis
there is such
\bdis
T_0[\vp_1;g]>0
\edis
that for every
\be \label{3.8}
T\in (T_0[\vp_1;g],+\infty)
\ee
the matrix
\be \label{3.9} \begin{split}
& \left\| \int_{\overset{p}{T}}^{\overset{s}{\wideparen{\overset{p-2}{T}+g}}}\prod_{r=0}^{s-1}
\tilde{Z}^2[\vp_1^r(t)]{\rm d}t\right\|_{p,s},\\
& p=s+1,\dots,k,\ s=1,\dots,k-1;\ 2\leq k\leq k_0
\end{split}
\ee
of the reversely iterated integrals (energies) fulfills the following equalities
\be \label{3.10}
\begin{split}
 & g=\int_{\overset{2}{T}}^{\overset{1}{\wideparen{\overset{1}{T}+g}}}\prod_{0}^{0} \\
 & \parallel \\
 & \int_{\overset{3}{T}}^{\overset{1}{\wideparen{\overset{2}{T}+g}}}\prod_{0}^{0} =
 \int_{\overset{3}{T}}^{\overset{2}{\wideparen{\overset{2}{T}+g}}}\prod_{0}^{1} \\
 & \parallel \\
 & \int_{\overset{4}{T}}^{\overset{1}{\wideparen{\overset{3}{T}+g}}}\prod_{0}^{0}=
 \int_{\overset{4}{T}}^{\overset{2}{\wideparen{\overset{2}{T}+g}}}\prod_{0}^{1}=
 \int_{\overset{4}{T}}^{\overset{3}{\wideparen{\overset{1}{T}+g}}}\prod_{0}^{2} \\
 & \parallel \\
 & \vdots \\
 & \parallel \\
 & \int_{\overset{k-1}{T}}^{\overset{1}{\wideparen{\overset{k-2}{T}+g}}}\prod_{0}^{0}=
 \int_{\overset{k-1}{T}}^{\overset{2}{\wideparen{\overset{k-3}{T}+g}}}\prod_{0}^{1}=
 \int_{\overset{k-1}{T}}^{\overset{3}{\wideparen{\overset{k-4}{T}+g}}}\prod_{0}^{2}=\dots =
 \int_{\overset{k-1}{T}}^{\overset{k-2}{\wideparen{\overset{1}{T}+g}}}\prod_{0}^{k-3} \\
 & \parallel \\
 & \int_{\overset{k}{T}}^{\overset{1}{\wideparen{\overset{k-1}{T}+g}}}\prod_{0}^{0}=
 \int_{\overset{k}{T}}^{\overset{2}{\wideparen{\overset{k-2}{T}+g}}}\prod_{0}^{1}=
 \int_{\overset{k}{T}}^{\overset{3}{\wideparen{\overset{k-3}{T}+g}}}\prod_{0}^{2}=\dots=
 \int_{\overset{k}{T}}^{\overset{k-2}{\wideparen{\overset{2}{T}+g}}}\prod_{0}^{k-3}=
 \int_{\overset{k}{T}}^{\overset{k-1}{\wideparen{\overset{1}{T}+g}}}\prod_{0}^{k-2}.
\end{split}
\ee
\end{mydef1}

\begin{remark}
Let us notice that the name \emph{local equalities} we use for every fixed row in the matrix (\ref{3.10}) only.
Namely, the corresponding segments
\bdis
[\overset{4}{T},\overset{1}{\wideparen{\overset{3}{T}+g}}],\
[\overset{4}{T},\overset{2}{\wideparen{\overset{2}{T}+g}}],\
[\overset{4}{T},\overset{3}{\wideparen{\overset{1}{T}+g}}],
\edis
for example, get on from the same point. Simultaneously, the distance of every two consecutive rows is (comp. (\ref{2.8}))
\bdis
\sim (1-c)\pi(T) \sim (1-c)\frac{T}{\ln T},\ T\to\infty.
\edis
\end{remark}

\begin{remark}
At the same time the name \emph{local equalities} is not quite exact. Namely, the elements of the matrix (\ref{3.10})
standing in the $2-(k-1)$th columns posses the property of internal non-locality (comp. part (b) of Remark 6).
\end{remark}

\section{Extension of non-local equalities (\ref{2.1}) by means of local equalities (\ref{3.10})}

Next, adding local equalities (\ref{3.10}) to non-local ones (\ref{2.1}) we obtain the following

\begin{mydef4}
\be \label{4.1}
\begin{split}
 & g=\int_{\overset{1}{T}}^{\overset{1}{\wideparen{T+g}}}\prod_{0}^{0} \\
 & \parallel \\
 & \int_{\overset{2}{T}}^{\overset{2}{\wideparen{T+g}}}\prod_{0}^{1}=
 \int_{\overset{2}{T}}^{\overset{1}{\wideparen{\overset{1}{T}+g}}}\prod_{0}^{0} \\
 & \parallel \\
 & \int_{\overset{3}{T}}^{\overset{3}{\wideparen{T+g}}}\prod_{0}^{2}=
 \int_{\overset{3}{T}}^{\overset{1}{\wideparen{\overset{2}{T}+g}}}\prod_{0}^{0}=
 \int_{\overset{3}{T}}^{\overset{2}{\wideparen{\overset{1}{T}+g}}}\prod_{0}^{1} \\
 & \parallel \\
 & \vdots \\
 & \parallel \\
 & \int_{\overset{k-1}{T}}^{\overset{k-1}{\wideparen{T+g}}}\prod_{0}^{k-2}=
 \int_{\overset{k-1}{T}}^{\overset{1}{\wideparen{\overset{k-2}{T}+g}}}\prod_{0}^{0}=
 \int_{\overset{k-1}{T}}^{\overset{2}{\wideparen{\overset{k-3}{T}+g}}}\prod_{0}^{1}= \dots =
 \int_{\overset{k-1}{T}}^{\overset{k-2}{\wideparen{\overset{1}{T}+g}}}\prod_{0}^{k-3} \\
 & \parallel \\
 & \int_{\overset{k}{T}}^{\overset{k}{\wideparen{T+g}}}\prod_{0}^{k-1}=
 \int_{\overset{k}{T}}^{\overset{1}{\wideparen{\overset{k-1}{T}+g}}}\prod_{0}^{0}=
 \int_{\overset{k}{T}}^{\overset{2}{\wideparen{\overset{k-2}{T}+g}}}\prod_{0}^{1}=\dots=
 \int_{\overset{k}{T}}^{\overset{k-2}{\wideparen{\overset{2}{T}+g}}}\prod_{0}^{k-3}=
 \int_{\overset{k}{T}}^{\overset{k-1}{\wideparen{\overset{1}{T}+g}}}\prod_{0}^{k-2}.
\end{split}
\ee
\end{mydef4}

\begin{remark}
Since (see (\ref{1.5}))
\be \label{4.2}
\begin{split}
 & \int_{\overset{p}{T}}^{\overset{s}{\wideparen{\overset{p-s}{T}+g}}}\prod_{r=0}^{s-1}\tilde{Z}^2
 [\vp_1^r(t)]{\rm d}t= \\
 & = \int_{\overset{p}{T}}^{\overset{s}{\wideparen{\overset{p-s}{T}+g}}}\prod_{r=0}^{s-1}
 \frac{\left|\zeta\left(\frac 12+i\vp_1^r(t)\right)\right|^2}
 {\omega[\vp_1^r(t)]}{\rm d}t=g
\end{split}
\ee
then we have the following: the set of equalities (\ref{4.1}) is a kind of set of constraints on
complicated behaviour of the function
\bdis
\zf,\ t\to\infty.
\edis
\end{remark}

\begin{remark}
Let us denote by
\bdis
S(T,g)
\edis
the complete set of equalities-constraints contained in the matrix (\ref{4.1}) for every fixed
\bdis
[T,g]:\ T\in (T_0[\vp_1,g],+\infty),\ g\in (0,\infty).
\edis
Since every integral (energy) from the matrix (\ref{4.1}) is invariant under the set of translations
\bdis
T\longrightarrow T';\ T,T'\in (T_0[\vp_1,g],+\infty)
\edis
then we call the set $S(T,g)$ \emph{the invariant set of equalities-constraints} for every
fixed $g\in (0,+\infty)$.
\end{remark}

\appendix

\section{}

In this part we give two examples on the properties of finite additivity and of finite multiplicativity
in the set of reversely iterated integrals (energies).

\subsection{On unbounded division of reversely iterated inegral (energy) on equal parts}

Since for every fixed
\bdis
[g,N]:\ g\in (0,+\infty),\ N\in\mbb{N};\ g=o\left(\frac{T}{\ln T}\right)
\edis
(see (\ref{1.1})) we have that
\bdis
g=N\delta(N)
\edis
then (see \cite{4}, (3.1), (3.3))
\be \label{a1}
\int_{\overset{k}{T}}^{\overset{k}{\wideparen{T+N\delta}}}\prod_{r=0}^{k-1}\tilde{Z}^2[\vp_1^r(t)]{\rm d}t=
N\int_{\overset{k}{T}}^{\overset{k}{\wideparen{T+\delta}}}\prod_{r=0}^{k-1}\tilde{Z}^2[\vp_1^r(t)]{\rm d}t .
\ee
Next, we have
\bdis
\begin{split}
 & \int_{\overset{k}{T}}^{\overset{k}{\wideparen{T+n\delta}}}=n\delta, \
 \int_{\overset{k}{T}}^{\overset{k}{\wideparen{T+(n-1)\delta}}}=(n-1)\delta \ \Rightarrow \
 \int_{\overset{k}{\wideparen{T+(n-1)\delta}}}^{\overset{k}{\wideparen{T+n\delta}}}=\delta, \\
 & n=1,\dots,N,
\end{split}
\edis
i.e.
\be \label{a2}
\int_{\overset{k}{\wideparen{T+(n-1)\delta}}}^{\overset{k}{\wideparen{T+n\delta}}}=
\int_{\overset{k}{T}}^{\overset{k}{\wideparen{T+\delta}}},\ n=1,\dots,N.
\ee

\begin{remark}
The formulae (\ref{a1}), (\ref{a2}) given a simultaneous $(k=1,\dots,k_0)$ unbounded division of the reversely
iterated integrals (energies) into equal parts.
\end{remark}

\subsection{A chain of integrals}

If we use finite case of the complete multiplicativity (see \cite{4}, (4.1)) together with the formulae
\bdis
\sqrt{\frac 2\pi}\int_0^\infty e^{-x^2/2}\cos\omega x{\rm d}x=e^{-\omega^2/2},\
\int_0^\infty e^{-x^2/2}{\rm d}x=\sqrt{\frac{\pi}{2}},
\edis
then we obtain the following result.

\begin{remark}
The following holds true
\bdis
\begin{split}
 & \int_{(\mbb{R}^+_0)^n}\left\{ \int_{\mbb{R}^+_0}\prod_{l=1}^n \cos(\omega_lx_l)\times \right. \\
 & \left. \times  \left[ \int_{\overset{k}{T}}^{\overset{k}{\wideparen{T+\prod\exp(-x_l^2/2)}}}
 \prod_{r=0}^{k-1}\tilde{Z}^2[\vp_1^r(t)]{\rm d}t\right]{\rm d}\omega_1\dots{\rm d}\omega_n\right\}
 {\rm d}x_1\dots{\rm d}\omega_n=\left(\frac{\pi}{2}\right)^n
\end{split}
\edis
for every fixed $n\in\mbb{N}$.
\end{remark}

\thanks{I would like to thank Michal Demetrian for his help with electronic version of this paper.}

\end{document}